%% file: root.tex
\documentclass{ifacconf}

\usepackage{graphicx}
\usepackage{natbib}        

\input{myheader}

\begin{document}
\begin{frontmatter}

\title{Null controllability of the 1D heat equation using flatness}


\author[First]{Philippe Martin}
\author[Second]{Lionel Rosier}
\author[First]{Pierre Rouchon}

\address[First]{Centre Automatique et Systèmes, MINES ParisTech,
75272 Paris, France (e-mail: \{philippe.martin,pierre.rouchon\}@mines-paristech.fr).}
\address[Second]{Institut Élie Cartan, UMR 7502 UdL/CNRS/INRIA, BP 70239,
54506 Vandœuvre-lès-Nancy, France (e-mail: Lionel.Rosier@univ-lorraine.fr)}

\begin{abstract}                
We derive in a straightforward way the null controllability of a 1-D heat equation with boundary control.
We use the so-called {\em flatness approach}, which consists in parameterizing the solution and the control by the derivatives of a ``flat output''. This provides an explicit control law achieving the exact steering to zero. We also give accurate error estimates when the various series involved are replaced by their partial sums, which is paramount for an actual numerical scheme. Numerical experiments  demonstrate the relevance of the approach.
\end{abstract}


\end{frontmatter}

\input{introduction}
\input{flatness}
\input{nullcont}
\input{estimates}
\input{numerics}

%


%
%
%

\bibliography{cpde}             

\end{document}

%% file: myheader.tex
\graphicspath{{Figures/}}

\usepackage[utf8]{inputenc}
\usepackage[T1]{fontenc}

\usepackage{subfigure}

\usepackage[retainorgcmds]{IEEEtrantools}
\usepackage{amssymb,amsmath}
\usepackage{eucal}

\providecommand{\abs}[1]{\left\lvert#1\right\rvert}
\providecommand{\norm}[1]{\left\lVert#1\right\rVert}

\newcommand{\R}{\mathbb R}
\newcommand{\C}{\mathbb C}
\newcommand{\N}{\mathbb N}

\newcommand{\be}{\begin{equation}}
\newcommand{\ee}{\end{equation}}

\newcommand{\bes}{\begin{equation*}}
\newcommand{\ees}{\end{equation*}}

\newcommand{\ba}{\begin{eqnarray}}
\newcommand{\ea}{\end{eqnarray}}
\renewcommand{\S}{\mathcal S}

\newcommand{\oi}{\overline{i}}

\newcommand{\ok}{\overline{k}}
\newcommand{\on}{\overline{n}}


%% file: introduction.tex
\section{Introduction}
%
The controllability of the heat equation was first considered in the 1-D case, \cite{FattoR1971ARMA,Jones1977JMAA,Littm1978ASNSPCS}), and very precise results were obtained by
the classical moment approach.  Next using Carleman estimates and duality arguments the null controllability was proved in \cite{FursiI1996book,LebeaR1995CPDE} for any bounded domain in~$\R^N$, any control time $T$, and any control region. This Carleman approach proves very efficient also with semilinear parabolic equations, \cite{FursiI1996book}.

By contrast the numerical control of the heat equation (or of parabolic equations) is in its early stage, see e.g.~\cite{MunchZ2010IP,BoyerHL2011NM,MicuZ2011SCL}. A natural candidate for the control input is the control of minimal $L^2-$norm, which may be obtained as a trace of the solution of the (backward) adjoint problem whose terminal state is the minimizer of a suitable quadratic cost. Unfortunately its computation is a hard task~\cite{MicuZ2011SCL}; indeed the terminal state of the adjoint problem associated with some regular initial state of the control problem may be highly irregular, which leads to severe troubles in the numerical computation of the control function.

All the above results rely on some observability inequalities for the adjoint system. A direct approach which does not involve the adjoint problem was proposed in~\cite{Jones1977JMAA,Littm1978ASNSPCS,LinL1995AMO}.  In \cite{Jones1977JMAA} a fundamental solution for the heat equation with compact support in time was introduced and used to prove null controllability. The results in \cite{Jones1977JMAA,Rosie2002CAM} can be used to derive control results on a bounded interval with one boundary control in some Gevrey class. An extension of those results to the semilinear heat equation in 1D was obtained in \cite{LinL1995AMO} in a more explicit way through the resolution of an ill-posed problem with data of Gevrey order 2 in~$t$.

In this paper we derive in a straightforward way the null controllability of the 1-D heat equation
\begin{IEEEeqnarray}{rCl'l}
    \theta_t(t,x) -\theta_{xx}(t,x) &=& 0, &(t,x)\in(0,T)\times(0,1) \label{B1}\\
    \theta_x(t,0) &=& 0, &t\in(0,T) \label{B2}\\
    \theta_x(t,1) &=& u(t), &t\in(0,T) \label{B3}
\end{IEEEeqnarray}
with initial condition
\begin{IEEEeqnarray*}{rCl'l}
    \theta(0,x) &=&\theta_0(x), &x\in (0,1). 
    \label{eq:ic}
\end{IEEEeqnarray*}
This system describes the dynamics of the temperature~$\theta$ in an insulated metal rod where the control~$u$
is the heat flux at one end. More precisely given any final time~$T>0$ and any initial state $\theta_0\in L^2(0,1)$
we provide an explicit control input $u\in L^2(0,T)$ such that the state reached at time~$T$ is zero, i.e.
\begin{IEEEeqnarray*}{rCl'l}
    \theta(T,x) &=&0, &x\in(0,1).
\end{IEEEeqnarray*}
We use the so-called {\em flatness approach}, \cite{FliesLMR1995IJoC}, which consists in parameterizing the solution~$\theta$ and the control~$u$ by the derivatives of a ``flat output''~$y$ (section~\ref{sec:flatness}); this notion was initially introduced for finite-dimensional (nonlinear) systems, and later extended to in particular parabolic PDEs, \cite{LarocMR2000IJRNC,LynchR2002IJC,MeureZ2008MCMDS,Meure2011A}. Choosing a suitable trajectory for this flat output~$y$ then yields an explicit series for a control achieving the exact steering to zero (section~\ref{sec:controllability}). This generalizes~\cite{LarocMR2000IJRNC}, where only approximate controllability was achieved through a similar construction. We then give accurate error estimates when the various series involved are replaced by their partial sums, which is paramount for an actual numerical scheme (section~\ref{sec:cestimates}). Numerical experiments  demonstrate the relevance of the approach (section~\ref{sec:numerics}).

In the sequel we will consider series with infinitely many derivatives of some functions. The notion of Gevrey order is a way of estimating the growth of these derivatives: we say that a function $y\in C^\infty([0,T])$ is {\em Gevrey of order $s\geq0$ on~$[0,T]$} if there exist positive constants~$M,R$ such that
\bes
\abs{y^{(p)}(t)} \leq M\frac{p!^s}{R^p} \qquad \forall t\in [0,T],\ \forall p\ge 0.
\ees
More generally if $K\subset\R^N$ is a compact set and $y$ is a function of class~$C^\infty$ on~$K$ (i.e. $y$ is the restriction to $K$ of a function of class~$C^\infty$ on some open neighbourhood $\Omega$ of~$K$), we say $y$ is {\em Gevrey of order $s_1$ in $x_1$, $s_2$ in $x_2$,\ldots,$s_N$ in $x_N$ on~$K$} if there exist positive constants $M,R_1,...,R_N$ such that
\bes
\abs{\partial_{x_1}^{p_1}\partial_{x_2}^{p_2}\cdots\partial_{x_N}^{p_N}y(x)} \le
M\frac{\prod_{i=1}^N (p_i !)^{s_i}}{\prod_{i=1}^N R_i^{p_i} }, \quad \forall x\in K,\ \forall p\in \N ^N.
\ees
By definition, a Gevrey function of order $s$ is also of order $r$ for~$r\geq s$. Gevrey functions of order~1 are analytic (entire if $s<1$).  Gevrey functions of order~$s>1$ have a divergent Taylor expansion; the larger~$s$, the ``more divergent'' the Taylor expansion. Important properties of analytic functions generalize to Gevrey functions of order $s>1$: the scaling, addition, multiplication and derivation of Gevrey functions of order $s>1$ is of order~$s$, see~\cite{Ramis1978,Rudin1987book}.  But contrary to analytic functions, functions of order $s>1$ may be constant on an open set without being constant everywhere.  For example the ``step function''
\bes
\phi_s(t):=\begin{cases}
1 & \text{if $t\leq0$}\\
0 & \text{if $t\geq1$}\\
\dfrac{ e^{-(1-t)^{-k}} }{ e^{-(1-t)^{-k}} + e^{-t^{-k}} }
&\text{if $t\in]0,1[$},
\end{cases}
\ees
where $k = (s-1)^{-1}$ is Gevrey of order~$s$ on $[0,1]$ (and in fact on~$\R$); notice $\phi_s(0)=1$, $\phi_s(1)=0$ and $\phi_s^{(i)}(0)=\phi_s^{(i)}(1)=0$ for all~$i\geq1$.

In conjunction with growth estimates we will repeatedly use Stirling's formula $n!\sim(n/e)^n\sqrt{2\pi n}$.

%% file: flatness.tex
\section{The heat equation is flat}\label{sec:flatness}
We claim the system~\eqref{B1}--\eqref{B3} is ``flat'' with $y(t):=\theta(0,t)$ as a flat output, which means there is (in appropriate spaces of smooth functions) a $1-1$ correspondence between arbitrary functions $t\mapsto y(t)$ and solutions of~\eqref{B1}--\eqref{B3}.

We first seek a formal solution in the form
\bes
\theta(t,x):=\sum_{i\ge0}\frac{x^i}{i!}a_i(t)
\ees
where the $a_i$'s are functions yet to define. Plugging this expression into~\eqref{B1} yields
\bes
\sum_{i\ge 0}\frac{x^i}{i!}[a_{i+2}-a_i']=0,
\ees
hence $a_{i+2}=a_i'$ for all~$i\ge0$. On the other hand $y(t)=\theta(0,t)=a_0(t)$, and \eqref{B2} implies~$a_1(t)=0$. As a consequence $a_{2i}=y^{(i)}$ and $a_{2i+1}=0$ for all~$i\ge0$. The formal solution thus reads
\be
\label{AA10}
\theta(t,x)=\sum_{i\ge 0}\frac{x^{2i}}{(2i)!}y^{(i)}(t)
\ee
while the formal control is given by
\be
\label{AA10bis}
u(t)=\theta_x(1,t)=\sum_{i\ge 1}\frac{y^{(i)}(t)}{(2i-1)!}.
\ee

We now  give a meaning to this formal solution by restricting $t\mapsto~y(t)$
to be Gevrey of order~$s\in [0,2)$.
\begin{prop}
\label{prop1}
Let $s\in [0,2)$, $-\infty <t_1<t_2<\infty$, and $y\in C^\infty ([t_1,t_2] )$ satisfying for some constants $M,R>0$
\be
\label{AA11}
\abs{y^{(i)}(t)} \le M \frac{i!^s}{R^i}, \qquad \forall i\ge 0,\ \forall t\in [t_1,t_2].
\ee
Then the function $\theta$ defined by \eqref{AA10} is Gevrey of order $s$ in $t$ and $s/2$ in $x$ on
$[t_1,t_2]\times [0,1]$; hence the control~$u$ defined by \eqref{AA10bis} is also Gevrey of order $s$ on $[t_1,t_2]$.
\end{prop}
\begin{pf}
We must prove the formal series
\be
\label{AA12}
\partial _t ^m\partial _x ^n \theta (t,x) = \sum_{2i\ge n}  \frac{ x^{2i-n } }{ (2i-n)! } y^{(i+n)}(t)
\ee
is uniformly convergent on $[t_1,t_2]\times [0,1]$ with growth estimates of the form
\be
\label{AA12bis}
\abs{\partial _t ^m\partial _x ^n \theta (t,x)} \le C
\frac{m! ^s}{R_1^m}\, \frac{n! ^\frac{s}{2}}{R_2^n} \cdot
\ee
By \eqref{AA11}, we have for all $(t,x)\in [t_1,t_2]\times [0,1]$
\begin{eqnarray*}
\left\vert  \frac{x^{2i-n}}{(2i-n)!} y^{(i+m)} (t) \right\vert &\le& \frac{M}{R^{i+m}} \, \frac{(i+m)! ^s}{(2i-n)!} \\
&\le& \frac{M}{R^{i+m} }\,  \frac{(2^{i+m} i! \, m! )^s}{(2i-n)!} \\
&\le& \frac{M}{R^{i+m}}\, \frac{2^{si} \bigl(2^{-2i} \sqrt{\pi i} \, (2i)! \bigr)^\frac{s}{2} }{(2i-n)!}\, \frac{m!^s}{2^{-sm}} \\
&\le&  M \frac{(\pi i) ^{\frac{s}{4}} }{R_1 ^i (2i-n)! ^{1-\frac{s}{2}} } n! ^{\frac{s}{2}} \frac{m!^s}{R_1 ^m},
\end{eqnarray*}
where we have set $R_1 = 2^{-s}R$; we have used Stirling's formula for $(2i)!$ and twice $(i+j)! \le 2^{i+j} i!j!$. Since $\sum_{2i\ge n} \frac{(\pi  i)^\frac{s}{4} }{R_1 ^i (2i-n)! ^{1-\frac{s}{2} }  } <\infty$ the series in~\eqref{AA12} are uniformly convergent for all $m,n\ge 0$, hence $\theta \in C^\infty ([t_1,t_2]\times [0,1])$.
Finally, since
\bes
\sum_{2i\ge n} \frac{ M(\pi  i)^\frac{s}{4} }{R_1 ^i (2i-n)! ^{1-\frac{s}{2} }  }
\le  M\Bigl(\frac{\pi}{2}\Bigr)^\frac{s}{4}   R_1 ^{-\frac{n}{2}}  \sum_{j\ge 0} \frac{ j^\frac{s}{4} + n^\frac{s}{4} }{R_1 ^\frac{j}{2}  j! ^{1-\frac{s}{2} }  }
\le C R_2^{-n}
\ees
where $R_2\in (0,\sqrt{R_1})$ and $C>0$ is some constant independent of~$n$, we have the desired estimates~\eqref{AA12bis}.\qed
\end{pf}

%% file: nullcont.tex
\section{Null controllability}\label{sec:controllability}
In this section we derive an explicit control steering the system from any initial state~$\theta_0\in L^2(0,1)$ at time~$0$ to the final state~$0$ at time~$T>0$. Two ideas are involved: on the one hand thanks to the flatness property it is easy to find a control achieving the steering to zero starting from a certain set of initial conditions (lemma~\ref{lem:steering}); on the other hand thanks to the regularizing property of the heat equation this set is reached from any~$\theta_0\in L^2(0,1)$ when applying first a zero control for some time (lemma~\ref{lem:zerocontrol}).

\begin{lem}\label{lem:steering}
Let $(y_i)_{i\ge0}$ be a sequence of real numbers such that for some constants~$M,R>0$
\be
\label{AA20}
\abs{y_i}\le M\frac{i!}{R^i}\qquad \forall i\ge 0.
\ee

Then the function defined on~$[t_1,t_1+R']$, $R'<R$, by
\bes
y(t):=\phi_s\Bigl(\frac{t-t_1}{R'}\Bigr)\sum_{i\ge0}y_i\frac{(t-t_1)^i}{i!},
\ees
is Gevrey of order~$s>1$ on $[t_1,t_1+R']$ and satisfies for all $i\geq0$
\begin{IEEEeqnarray}{rCl}
y^{(i)}(t_1)  &=& y_i\label{eq:yt1}\\
y^{(i)}(t_1+R')  &=& 0\label{eq:yt2}.
\end{IEEEeqnarray}

Moreover the control defined on~$[t_1,t_1+R']$ by
\be\label{AA31}
u(t):=\sum_{i\ge1}\frac{y^{(i)}(t)}{(2i-1)!}
\ee
is also Gevrey of order $s$ on $[t_1,t_1+R']$ and steers the system from the initial state $\sum_{i\ge0}y_i\frac{x^{2i}}{(2i)!}$ at time~$t_1$ to the final state~$0$ at time~$t_1+R'$.
\end{lem}



\begin{pf}
Let $\overline{y}(z):=\sum_{i\geq0}y_i\frac{z^i}{i!}$. The growth property~\eqref{AA20} implies $\overline{y}$ is analytic on the disc $\{z\in\C;\ \abs z<R\}$, the convergence being moreover uniform for $\abs z\leq R'<R$. Therefore $\overline{y}$ is Gevrey of order~$1$, hence of order~$s>1$, on~$[t_1,t_1+R']$. On the other hand $\phi(t):=\phi_s\Bigl(\frac{t-t_1}{R'}\Bigr)$ is also Gevrey of order~$s$ on $[t_1,t_1+R']$, hence so is the product~$y$ of $\overline y$ and~$\phi$. The boundary values~\eqref{eq:yt1}-\eqref{eq:yt2} follow at once from the definition of~$\phi_s$.

The control~$u$ in~\eqref{AA31} achieves the steering to zero; indeed
\bes
\theta(t,x):=\sum_{i\ge 0}\frac{ x^{2i}}{(2i)!}y^{(i)}(t),
\ees
as well as~$u$, is by proposition~\ref{prop1} Gevrey of order~$s$ in~$t$ and $s/2$ in~$x$, and obviously satisfies~$\theta(t_1+R',x)=0$. \qed
\end{pf}

\begin{lem}\label{lem:zerocontrol}
Let $\theta_0\in L^2(0,1)$ and $\tau>0$. Consider the final state $\theta_\tau(x):=\theta(\tau,x)$ reached when applying the control $u(t):=0$, $t\in[0,\tau]$, starting from the initial state~$\theta_0$.

Then $\theta_\tau$ is analytic in $\C$ and can be expanded as
\bes
\theta_\tau(x) = \sum_{i\ge0}y_i\frac{x^{2i}}{(2i)!}, \quad x\in\C, 
\ees
with
\bes
\abs{y_i} \le C\Bigl(1+\frac{1}{\sqrt\tau}\Bigr)\frac{i!}{\tau^i}
\ees
where $C$ is some positive constant depending only on~$\theta_0$.
\end{lem}

\begin{pf}
Decompose $\theta _0$ as the Fourier series of cosines
\bes
\theta_0(x)=\sum_{n\ge 0}c_n\sqrt{2}\cos(n\pi x)
\ees
where the convergence holds in $L^2(0,1)$ and
\bes
2|c_0|^2+\sum_{n\ge 1} |c_n|^2 =\int_0^1 |\theta _0(x)|^2 dx <\infty.
\ees
The solution starting from~$\theta_0$ then reads
\be
\theta (t,x)=\sum_{n\ge 0}c_ne^{-n^2\pi^2t}\sqrt{2}\cos(n\pi x)
\label{C100}
\ee
and in particular
\bes
\theta_\tau(x)=\sum_{n\ge 0}c_n e^{-n^2\pi^2\tau}\sqrt{2}\cos(n\pi x).
\ees

The series for $\theta_\tau$ is analytic in~$\C$ since for all~$\abs{x}\leq r$
\bes
\abs{c_ne^{-n^2\pi^2\tau}\sqrt{2}\cos(n\pi x)} \le C_1\left(\sup_{k\ge0}\abs{c_k}\right)e^{-n^2\pi^2\tau+n\pi r}
\ees
where $C_1$ is some positive constant; this ensures the uniform convergence of the series in every open disk of radius $r>0$.

Moreover
\begin{IEEEeqnarray*}{rCl}
\theta_\tau(x) &=& \sqrt{2}\sum_{n\ge 0}c_ne^{-n^2\pi^2\tau}\sum_{i\ge0}(-1)^i\frac{(n\pi x)^{2i}}{(2i)!}\\
%
%
&=& \sum_{i\ge0} \frac{x^{2i}}{(2i)!} \underbrace{\left(\sqrt{2}(-1)^i\sum_{n\ge0}c_n e^{-n^2\pi^2\tau}(n\pi)^{2i}\right)}_{=:y_i}
\end{IEEEeqnarray*}
The change in the order of summation will be justified once we have proved that $y_i$, $i\geq0$, is absolutely convergent and
\bes
\sum_{i\ge0}\abs{y_i}\frac{x^{2i}}{(2i)!} <\infty,  \quad\forall x\geq0.
\ees
For $i\ge 0$ let $h_i(x):=e^{-\tau\pi^2x^2}(\pi x)^{2i}$ and $N_i:=\left[\bigl(\frac{i}{\pi ^2\tau}\bigr)^{\frac{1}{2}}\right]$. The map $h_i$ is increasing on $\bigl[0,\left(\frac{i}{\pi ^2\tau}\right)^{\frac{1}{2}}\bigr]$ and decreasing on $\bigl[\left(\frac{i}{\pi ^2\tau}\right)^{\frac{1}{2}},+\infty\bigr)$ hence
\begin{IEEEeqnarray*}{rCl}
\sum_{n\ge0}h_i(n) &\le& \int_0^{N_i}h_i(x)dx + h_i(N_i)\\
&& +\: h_i (N_i + 1) + \int_{N_{i+1} }^\infty h_i(x)dx\\
&\le& 2h_i\left(\Bigl(\frac{i}{\pi^2\tau}\Bigr)^{\frac{1}{2}}\right) + \int_0^\infty h_i(x)dx\\
&\le& C_2\frac{i!}{\tau^i\sqrt{i}} + \int_0^\infty h_i(x)dx;
\end{IEEEeqnarray*}
$C_2$ is some positive constant and we have used Stirling's formula. On the other hand integrating by parts yields
\begin{IEEEeqnarray*}{rCl}
\int_0^\infty h_i(x)dx &=& \frac{2i-1}{2\tau}\int_0^\infty h_{i-1}(x)dx\\
&=& \frac{(2i-1)\cdots 3\cdot1}{(2\tau)^i}\int_0^\infty e^{-\tau\pi^2x^2}dx\\
&=& \frac{(2i)!}{2^ii!(2\tau)^i}\cdot\frac{1}{\pi\sqrt\tau} \int_0^\infty e^{-x^2}dx\\
&\leq& C_3\frac{i!}{\tau^i\sqrt{i\tau}},
\end{IEEEeqnarray*}
where $C_3$ is some positive constant and we have again used Stirling's formula. As a consequence
\be\label{Z7}
\abs{y_i}\le\sqrt{2}\sup_{n\ge 0}\abs{c_n}\sum_{n\ge0}h_i(n) \le C\Bigl(1+\frac{1}{\sqrt\tau}\Bigr)\frac{i!}{\tau^i}
\ee
where $C$ is some positive constant. Finally
\bes
\sum_{i\ge0}\abs{y_i}\frac{x^{2i}}{(2i)!}
\leq C\Bigl(1+\frac{1}{\sqrt\tau}\Bigr)\sum_{i\ge0}\underbrace{\frac{i!}{(2i)!}\Bigl(\frac{x^2}{\tau}\Bigr)^i}_{=:v_i}<\infty
\ees
since $\frac{v_{i+1}}{v_i}\sim\frac{1}{4i}\frac{x^2}{\tau}$.\qed
\end{pf}

With the two previous lemma at hand we can now state our main controllability result.
\begin{thm}\label{thm1}
Let $\theta_0\in L^2(0,1)$ and~$T>0$. Pick any $\tau\in(0,T)$ and $s\in (1,2)$. Then there exists a function~$y$
Gevrey of order~$s$ on~$[\tau,T]$ such that the control
\bes
u(t):=\begin{cases}
0 & \text{if $0\le t\le\tau$}\\
\sum_{i\ge1}\frac{y^{(i)}(t)}{(2i-1)!} &\text{if $\tau <t\le T$}.
\end{cases}
\ees
steers the system from the initial state $\theta_0$ at time~$0$ to the final state~$0$ at time~$T$.

Moreover $u$ is Gevrey of order~$s$ on~$[0,T]$; $t\mapsto\theta(t,\cdot)$ is in $C\bigl([0,T],L^2(0,1)\bigr)$; $\theta$ is Gevrey of order~$s$ in~$t$ and~$s/2$ in~$x$ on~$[\varepsilon,T]\times[0,1]$ for all $\varepsilon\in(0,T)$.
\end{thm}
\begin{pf}
By lemma~\ref{lem:zerocontrol} the state reached at time~$\tau$ reads $\sum_{i\ge0}y_i\frac{x^{2i}}{(2i)!}$ with the sequence
$(y_i)_{i\ge0}$ satisfying the growth property of lemma~\ref{lem:steering} with $M:=1+\frac{1}{\sqrt\tau}$ and~$R:=\tau$. Hence the desired function is given by
\bes
y(t):=\phi_s\Bigl(\frac{t-\tau}{R'}\Bigr)\sum_{i\ge0}y_i\frac{(t-\tau)^i}{i!},
\ees
on~$[\tau,\tau+R']$, where $R'<\tau$ and $R'\leq T-\tau$; and by $y(t):=0$ on~$[\tau+R',T]$. Moreover $y$ and~$u$ are Gevrey of order~$s$ on~$[\tau,T]$, and by construction $\theta(\tau,x)=\theta(\tau^+,x)=\sum_{i\ge0}y_i\frac{x^{2i}}{(2i)!}$ and $u(\tau)=u(\tau^+)=0$.
By Proposition~\ref{prop1} the solution~$\theta$ on~$[\tau,T]\times[0,1]$ is well-defined and Gevrey of order~$s$ in~$t$ and $s/2$ in~$x$ hence $t\mapsto\theta(t,\cdot)\in C((0,T];C^1([0,1]))$ and $u\in C([0,T])$.
On the other hand it is easily seen that $\theta$ is Gevrey of order~$1$ in~$t$ and $1/2$ in~$x$ on~$[\varepsilon,\tau]\times[0,1]$ for all $\varepsilon\in (0,\tau )$.
%

Thus the solution $\theta$ is Gevrey of order $1$ in $t$ and $1/2$ in $x$ on $[\varepsilon , \tau ]\times [0,1]$ while it is  Gevrey of order $s$ in $t$ and $s/2$ in $x$ on $[\tau ,T]\times [0,1]$.
To prove $\theta$ is Gevrey of order $s$ in $t$ and $s/2$ in $x$ on $[\varepsilon , T]\times [0,1]$ it is then sufficient to check
$\partial _t^k \theta (\tau, x ) = \partial _t ^k \theta (\tau ^+, x)$ for $k\geq0$ and~$x\in [0,1]$. But
\begin{IEEEeqnarray*}{rCl}
\partial _t^k \theta (\tau ^+, x ) &=& \sum_{i\ge 0} \frac{x^{2i}}{(2i)!} y^{(i+k)} (\tau ) \\
&=& \sum_{i\ge 0} \frac{x^{2i}}{(2i)!} y_{i+k} \\
&=& \sqrt{2} \sum_{i\ge 0} \frac{x^{2i}}{(2i)!}
\left( \sum_{n\ge 0} c_n e^{-n^2\pi ^2 \tau } n^{2(i+k)}  \right) (-\pi^2)^{i+k} \\
&=& \sum_{n\ge 0}  c_n (-n^2\pi ^2) ^k e^{-n^2\pi ^2 \tau }  \sqrt{2} \cos (n\pi x) \\
&=& \partial _t ^k \theta (\tau, x).
\end{IEEEeqnarray*}
As a consequence $u$ is also Gevrey of order~$s$ on~$[0,T]$. \qed
\end{pf}

%% file: estimates.tex
\section{Numerical estimates}\label{sec:cestimates}
Summarizing the previous section the control~$u$ and solution~$\theta$ on~$[\tau,\tau+R']$ are given by the infinite series
\begin{IEEEeqnarray}{rCl}
\label{eq:u}u(t) &=& \sum_{i\ge1}\frac{y^{(i)}(t)}{(2i-1)!}\\
\theta(t,x) &=& \sum_{i\ge 0}y^{(i)}(t)\frac{x^{2i} }{(2i)!}\\
y(t) &=& \phi_s\Bigl(\frac{t-\tau}{R'}\Bigr)\sum_{k\ge0}y_k\frac{(t-\tau)^k}{k!}\\
\label{eq:y}y_k &=& \sqrt{2}\Bigl(\sum_{n\ge 0} c_{n} e^{-n^2\pi ^2 \tau } n^{2k} \Bigr) (-\pi^2)^k;
\end{IEEEeqnarray}
moreover $y$ hence $u$ and~$\theta$ are identically zero on~$[\tau+R',T]$.
The aim of this section is to show that the partial sums
\begin{IEEEeqnarray}{rCl}
\label{eq:ubar}\overline{u}(t) &:=& \sum_{1\le i\le \oi}\frac{y^{(i)}(t)}{(2i-1)!}\\
\label{eq:thetabar}\overline{\theta}(t,x) &:=& \sum_{0\le i\le \oi}y^{(i)}(t)\frac{x^{2i} }{(2i)!}\\
\label{eq:ybar}\overline{y}(t) &:=& \phi_s\Bigl(\frac{t-\tau}{R'}\Bigr)\sum_{0\le k\le\ok}y_k\frac{(t-\tau)^k}{k!}\\
\label{eq:ykbar}\overline{y_k} &:=& \sqrt{2}\Bigl(\sum_{0\le n\le \on} c_{n} e^{-n^2\pi ^2 \tau } n^{2k} \Bigr) (-\pi^2)^k.
\end{IEEEeqnarray}
for given $\oi,\ok,\on\in\N$ provide very good approximations of the above series, and to give explicit error estimates.

\begin{thm}
\label{thm4}
There exist positive constants $C,C_1,C_2,C_3$ such that for all
$\theta_0\in L^2$, $\oi,\ok,\on\in\N$, and~$t\in[\tau,T]$
\bes
\norm{\theta(t)-\overline{\theta} (t)}_{L^\infty}
\le C\left( e^{-C_1 \, \oi \ln \oi} + e^{-C_2\,  \ok} + e^{-C_3\, \on ^2} \right)
\norm{\theta_0}_{L^2}
\ees
\end{thm}
\begin{pf}
First notice that for $(t,x)\in [\tau,T]\times[0,1]$
\bes
\label{delta}
\abs{\theta (t,x) - \overline{\theta} (t,x)} \le \Delta_1 + \Delta_2 + \Delta_3,
\ees
where
\begin{IEEEeqnarray*}{rCl}
\Delta_1 &:=& \Bigl\vert \sum_{i>\oi}y^{(i)}(t)\frac{x^{2i}}{(2i)!} \Bigr\vert\\
\Delta_2 &:=& \Bigl\vert \sum_{0\le i\le \oi } \partial_t ^i \Big[ \phi (t)
\sum_{k > \ok } y_{k}  \frac{(t-\tau )^k}{k!} \Big] \frac{x^{2i}}{(2i)!} \Bigr\vert\\
\Delta_3 &:=& \Bigl\vert \sum_{0\le i\le \oi } \partial _t ^i \Big[ \phi (t)
\sum_{0\le k \le \ok }
\\&&\quad
\sqrt{2} \big(\sum_{n> \on} c_n e^{-n^2\pi ^2 \tau } n^{2k}\big) (-\pi^2)^k \frac{(t-\tau )^k}{k!} \Big] \frac{x^{2i}}{(2i)!} \Bigr\vert .
\end{IEEEeqnarray*}

By lemma~\ref{lem:steering} $y$ is Gevrey of order~$s$ on~$[\tau,T]$ with some~$M_1\norm{\theta_0}_{L^2},R_1>0$ hence
\begin{IEEEeqnarray*}{rCl}
\Delta_1 &\leq& \sum_{i>\oi}\frac{\abs{y^{(i)}(t)}}{(2i)!}\\
&\leq& M_1\norm{\theta_0}_{L^2}\sum_{i>\oi}\frac{i!^s}{(2i)!R_1^i}\\
&\leq& M_1'\norm{\theta_0}_{L^2}\sum_{i>\oi}\frac{\sqrt i}{(4R_1)^i\sqrt i^{2-s}} \Bigl(\frac{i}{e}\Bigr)^{(s-2)i}\\
&\leq& M_1'\norm{\theta_0}_{L^2}\sum_{i>\oi}\frac{ i^{\frac{s-1}{2}}e^{(2-s)(1-\ln i)i} }{(4R_1)^i},
\end{IEEEeqnarray*}
where we have used Stirling's formula. Pick $C_1<2-s$ and $\sigma\in (C_1, 2-s)$. Then for $i > \oi$
\bes
(4R_1)^{-i} i^{\frac{s-1}{2}} e^{(2-s) (1-\ln i )i } \le K_1e^{-\sigma i (\ln i -1 )},
\ees
where $K_1=K_1(s, \sigma , R_1)$. But
\begin{IEEEeqnarray*}{rCl}
\sum_{i > \oi} e^{-\sigma i(\ln i -1 ) } &\le& \int_{\oi} ^\infty e ^{-\sigma x (\ln x -1 ) } dx\\
&\leq& K_1' \int _{\oi (\ln \oi -1 )} ^\infty e ^{-\sigma x} dx\\
&\leq& K_1' e^{-C_1 \oi \ln \oi},
\end{IEEEeqnarray*}
so that eventually
\bes
\Delta_1 \le K_1''\norm{\theta_0}_{L^2} e^{-C_1\oi\ln \oi}.
\ees

For $\Delta_2$ we first notice that for $t\in\{ z\in \C; |z-\tau |\le \rho'\tau \}$, where~$\rho'$ satisfies
$R'/\tau<\rho'<1$,
\begin{IEEEeqnarray*}{rCl}
\big\vert \sum_{k>\ok} y_k \frac{( z -\tau ) ^k}{k!} \big\vert
&\le& \sum_{k>\ok} |y_k| \frac{{\rho '} ^k \tau ^k }{k!}\\
&\leq& M_1\norm{\theta_0}_{L^2}\sum_{k>\ok}{\rho'}^k\\
&\leq& M_1''\norm{\theta_0}_{L^2}{\rho'}^{\ok+1}.
\end{IEEEeqnarray*}
By the Cauchy estimates we thus have for $\tau\le t\le\tau+R'$
\bes
\Big\vert  \partial_t^i\big[\sum_{k> \ok } y_{j,k} \frac{(t-\tau ) ^k }{k!}\big]  \Big\vert
\le M_1''\norm{\theta_0}_{L^2}{\rho'}^{\ok + 1}  \frac{i!}{R_1^i},
\ees
hence
\bes
\Big\vert\partial _t ^i \big[\phi (t) \sum_{k>\ok} y_k\frac{(t-\tau )^k}{k!}\big] \Big\vert
\le M_2\norm{\theta_0}_{L^2}{\rho'}^{\ok + 1} \frac{i! ^s}{R_2^i}.
\ees
It follows that
\bes
\Delta_2
\le \sum_{0\le i\le \oi} M_2\norm{\theta_0}_{L^2}{\rho '}^{\ok + 1}
\frac{i!^s}{(2i)! R_3 ^i}
\le K_2\norm{\theta_0}_{L^2}e^{-C_2\ok}
\ees
where $0<C_2 <\ln\frac{\tau}{R'}$.

To estimate $\Delta_3$ first notice that for~$\alpha>0$
\bes
\sum_{n>\on}e^{-\alpha n^2} \leq \int_{\on}^\infty e^{-\alpha x^2}dx
= \int_{\on^2}^\infty\frac{e^{-\alpha y}}{2\sqrt y}dy
\leq
\frac{e^{-\alpha n^2}}{2\alpha\on}.
\ees

Pick $\rho '$ and $\rho ''$ with $R\/\tau < \rho ' <\rho '' <1$. Then
\begin{IEEEeqnarray*}{rCl}
\Big\vert\sum_{n>\on}c_n e^{-n^2\pi ^2 \tau } n^{2k}\Big\vert
&\leq& \norm{\theta_0}_{L^2}\frac{e^{-k}k^k}{(\pi^2\rho''\tau)^k}
\sum_{n>\on}e^{-n^2\pi^2(1-\rho'')\tau}\\
&\leq& K_3\norm{\theta_0}_{L^2}\frac{k!}{(\pi^2\rho''\tau)^k}
\frac{e^{-\pi^2(1-\rho'')\tau\on^2}}{\on},
\end{IEEEeqnarray*}
where we have used first that $x\to e^{-x}x^k$ is maximum for $x=k$, and then that $e^{-k}k^k\le Ck!$ by Stirling's formula.
Hence for $z\in\{ z\in\C;\ \ |z- \tau | \le \rho'\tau \}$
\begin{IEEEeqnarray*}{rCl}
\IEEEeqnarraymulticol{3}{l}{
\Big\vert\sum_{0\le k\le \ok}\sqrt{2}
\Big( \sum_{n>\on} c_{j,n} e^{-n^2\pi ^2 \tau} n^{2k}\Big)  (-\pi^2)^k\frac{( z -\tau )^k}{k!}  \Big\vert}\\ 
\qquad
&\le& K_3'\norm{\theta_0}_{L^2}\sum_{0\le k\le\ok} 
\frac{e^{-\pi^2(1-\rho'')\tau\on^2}}{\on}  \Big(\frac{\rho'}{\rho''}\Big)^k\\
&\le& K_3''\norm{\theta_0}_{L^2}\frac{e^{-\pi^2(1-\rho'')\tau\on^2}}{\on}.
\end{IEEEeqnarray*}

Setting $C_3:=\pi^2(1-\rho'')\tau<\pi^2(\tau-R')$, this yields for $\tau\le t\le \tau+R'$
\begin{IEEEeqnarray*}{rCl}
\IEEEeqnarraymulticol{3}{l}{
\Big\vert
\partial _t^i \big[ \sum_{0\le k\le \ok} \sqrt{2} \big( \sum_{n> \on} c_n e^{-n^2\pi ^2 \tau} n^{2k}  \big)  (-\pi^2)^k\frac{(t-\tau )^k}{k!}\big]
\Big\vert}\\ 
\qquad
&\le& M_3\norm{\theta_0}_{L^2} e^{-C_3\on ^2 } \frac{i!}{R_1^i}
\end{IEEEeqnarray*}
and
\begin{IEEEeqnarray*}{rCl}
\IEEEeqnarraymulticol{3}{l}{
\Big\vert
\partial _t^i \big[\phi (t)  \sum_{0\le k\le \ok} \sqrt{2} \big( \sum_{n> \on} c_{j,n} e^{-n^2\pi ^2 \tau} n^{2k}  \big)  (-\pi^2)^k\frac{(t-\tau )^k}{k!} \big]
\Big\vert}\\
\qquad
&\le& M_3'\norm{\theta_0}_{L^2} e^{-C_3\on ^2 } \frac{i! ^s}{R_2^i}.
\end{IEEEeqnarray*}

We then conclude
\bes
\Delta_3 \le \sum_{0\le i\le \oi} M_3'\norm{\theta_0}_{L^2} e^{-C_3\on^2} \frac{i!^s}{(2i)! R_2^i}.
\le M_3''\norm{\theta_0}_{L^2}e^{-C_3\on ^2}
\label{delta4}
\ees

Collecting the inequalities for $\Delta_1,\Delta_2,\Delta_3$ eventually gives the statement of the theorem.\qed
\end{pf}

Let $\hat\theta$ denote the solution of \eqref{B1}--\eqref{B3} where the ``exact'' control~$u$~\eqref{eq:u} is replaced by the truncated control~$\overline{u}$~\eqref{eq:ubar}, still starting from the initial condition~$\theta_0$. Notice $\hat\theta$ is the ``exact'' solution obtained when applying the truncated control, while $\overline{\theta}$ is the truncated solution obtained when applying the ``exact'' control. Since in ``real life'' only the truncated control can be actually applied it is important to know how well $\hat\theta$ approximates~$\theta$, the ``exact'' solution when applying the ``exact'' control. It turns that $\hat\theta$ satisfies the same relation as~$\overline\theta$ in Theorem~\ref{thm4}; the proof is omitted for lack of space but follows from the proof of Theorem~\ref{thm4}.
\begin{cor}\label{cor3}
With the same notations as in Theorem \ref{thm4},
%
\bes
\norm{\theta(t)-\hat{\theta} (t)}_{L^\infty}
\le C\left( e^{-C_1 \, \oi \ln \oi} + e^{-C_2\,  \ok} + e^{-C_3\, \on ^2} \right)
\norm{\theta_0}_{L^2}
\ees
\end{cor}

%% file: numerics.tex
\section{Numerical experiments}\label{sec:numerics}
\begin{figure}[ht!]
\centering
\includegraphics[width=1.0\columnwidth]{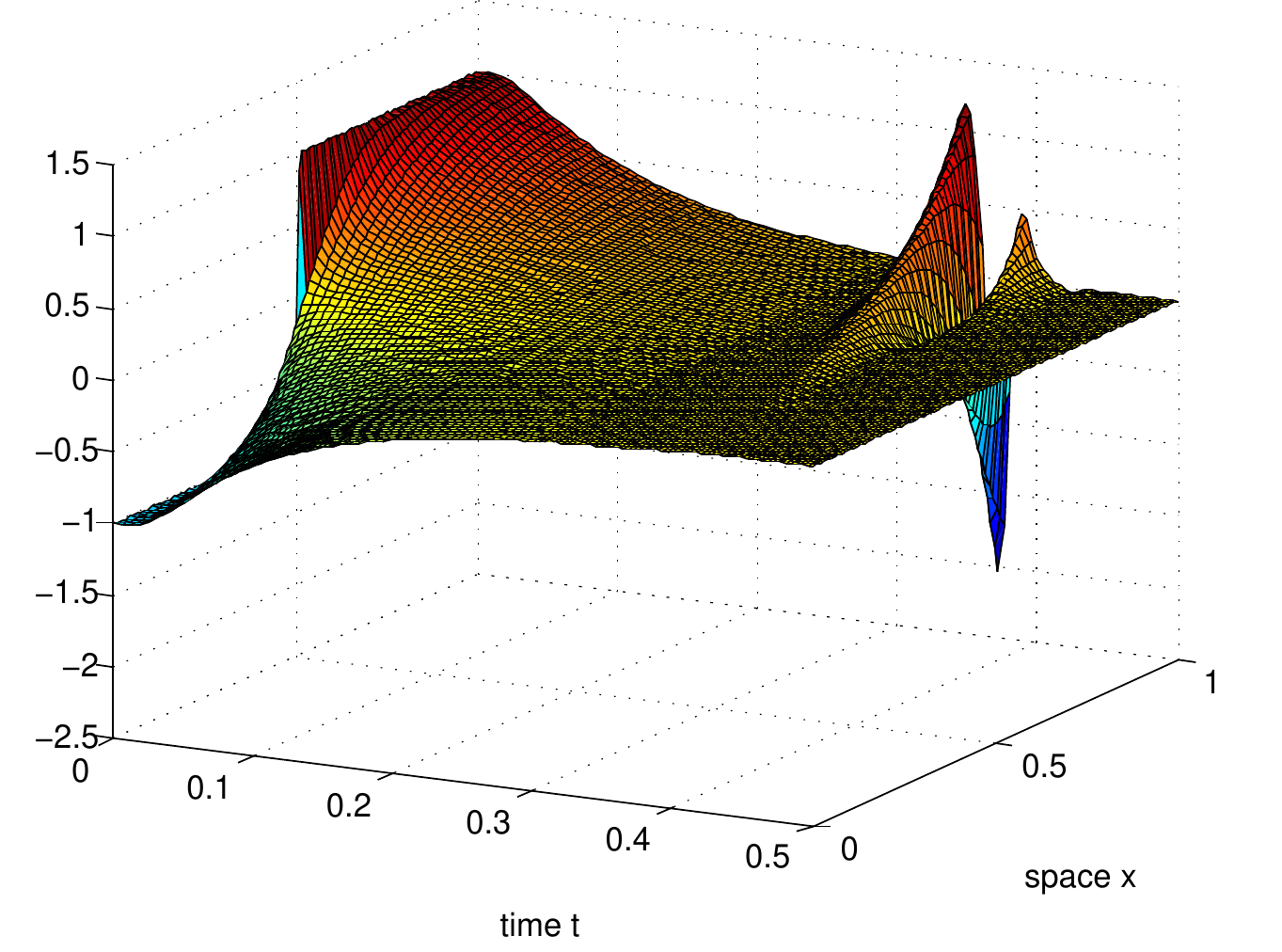}
\caption{$\overline\theta(t,x)$ for $R'=0.2$ and $s=1.6$.}
\label{fig:surf1130}
\end{figure}

\begin{figure}[ht!]
\centering
\includegraphics[width=1.0\columnwidth]{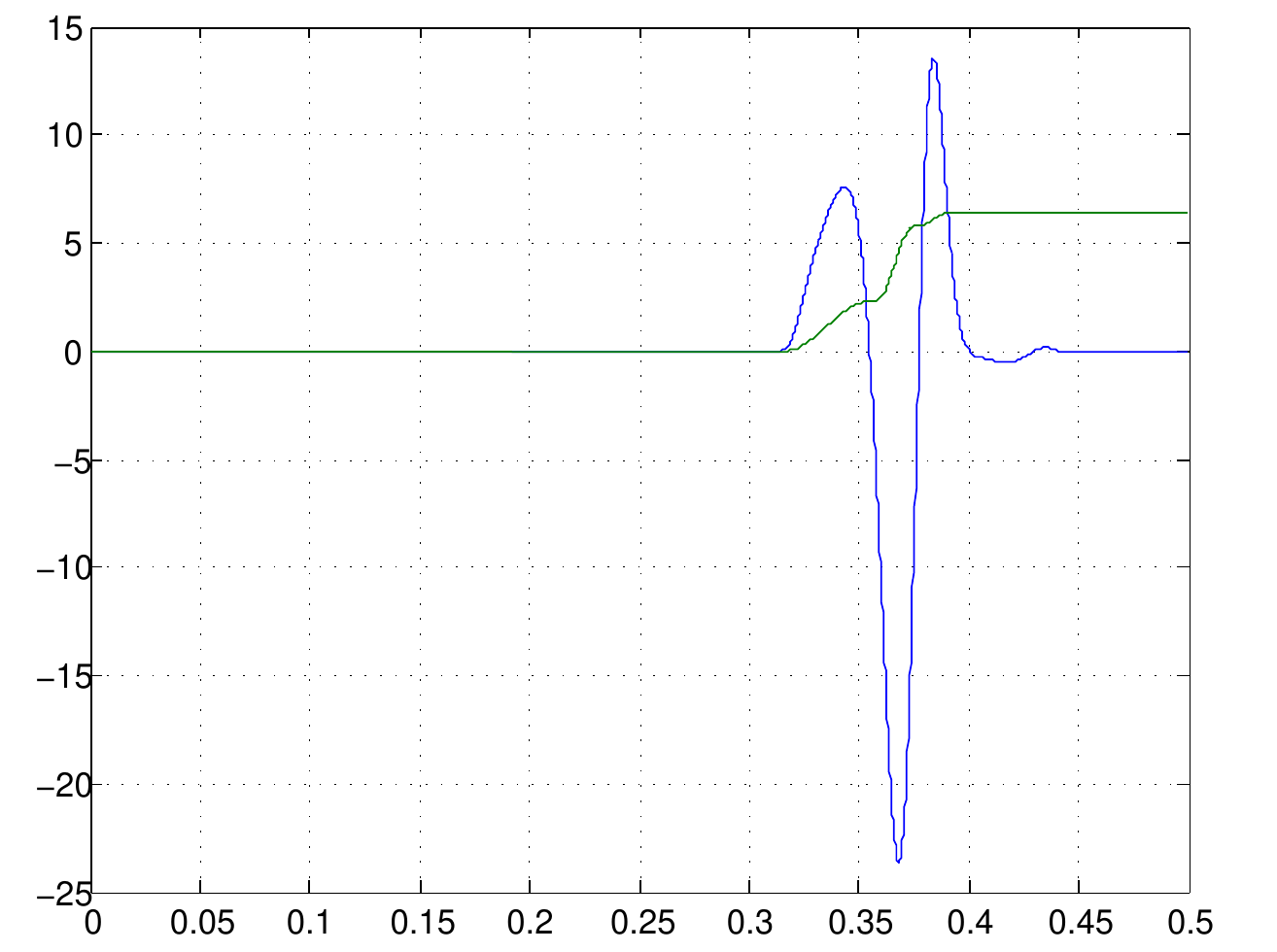}
\caption{$\overline u(t)$ and $\norm{\overline u}_{L^2(0,t)}$ for $R'=0.2$ and $s=1.6$.}
\label{fig:control1130}
\end{figure}

We hare conducted a number of numerical experiments to demonstrate the relevance of our approach, and to investigate the influence of the main parameters. We have focused on the control effort, which is probably the most important quantity. The results are summarized in the two tables below, where $\norm{\overline u}_{L^2}$ (top table) and $\norm{\overline u}_{L^\infty}$ (bottom table) are given for different values of the parameters $s$ and~$R'$ in~\eqref{eq:ybar}. Figures~\ref{fig:surf1130} and~\ref{fig:control1130} give moreover the complete temperature and control evolution for the case $R'=0.2$ and~$s=1.6$. For all experiments the ``regularization time'' $\tau$ is~$0.3$, and the initial condition $\theta_0$ is a step function with $\theta_0(x)=-1$ on~$[0,1/2)$ and $\theta_0(x)=1$ on~$(1/2,1)$ hence Fourier coefficients $c_{2p}=0$ and $c_{2p+1}=\frac{(-1)^{p+1}}{2p+1}\frac{2\sqrt2}{\pi}$ for $p\geq0$; $\oi,\ok,\on$ are ``large enough'' for a good accuracy of the truncated series\eqref{eq:ubar}~--\eqref{eq:ykbar}.

\begin{center}
\begin{tabular}{|r|c|c|c|c|c|}\hline
$\mathbf{s,R'}$& \textbf{0.15}& \textbf{0.20} & \textbf{0.25} & \textbf{0.30}\\ \hline
\textbf{1.5}& 693& 63.3& 12.7& 3.82\\
\hline
\textbf{1.6}& 35.3& 6.41& 1.95& 0.78\\
\hline
\textbf{1.7}& 7.49& 1.95& 0.74& 0.34\\
\hline
\textbf{1.8}& 5.53& 1.24& 0.48& 0.23\\
\hline
\textbf{1.9}& 5.71& 1.29& 0.47& 0.22\\
\hline
\end{tabular}
\end{center}

\begin{center}
\begin{tabular}{|r|c|c|c|c|c|}\hline
$\mathbf{s,R'}$& \textbf{0.15}& \textbf{0.20} & \textbf{0.25} & \textbf{0.30}\\ \hline
\textbf{1.5}& 3666& 330& 55.2& 18.1\\
\hline
\textbf{1.6}& 118& 23.6& 7.17& 2.76\\
\hline
\textbf{1.7}& 18.6& 4.78& 1.73& 0.76\\
\hline
\textbf{1.8}& 36.4& 4.59& 1.44& 0.65\\
\hline
\textbf{1.9}& 47.8& 9.66& 2.13& 0.82\\
\hline
\end{tabular}
\end{center}

A few conclusions can be drawn from these experiments:
\begin{itemize}
  \item when $R'$ gets small there seems to be an ``optimal'' value of~$s$ which minimizes the control effort (the values are different for $\norm{\overline u}_{L^2}$ and $\norm{\overline u}_{L^\infty}$)
  \item when $s$ it too small or too large the derivatives $\overline y^{(i)}$ tend to ``crowd'' near the extremities of~$[\tau,\tau+R']$ with large amplitudes
  \item the smaller $R'$ the more noticeable this phenomenon, with of course an increased control effort.
\end{itemize}

An interesting question is the tradeoff between $\tau$ and~$R'$ to reach the zero state at time~$T:=\tau+R'$ with the smallest control effort: longer regularization~$\tau$ or longer duration~$R'$ of the active control? With the present construction of the function~$y$~\eqref{eq:y} we have by design $R'\leq\tau$, which is probably an important restriction. Other, but more complicated, constructions without this limitation are possible and will be studied in the future.

